# Hasse theorem__ an elementary approach

Jianhua Chen, Debiao He *, Zhijin Hu, Yitao Chen, Hao Hu

School of mathematics and statistics Wuhan university, Wuhan, China

*Email: hedebiao@163.com

Abstract. We give an elementary proof to Hasse theorem.

Let $k$ be arbitrary field, possibly of characteristic $p \neq 2, 3$, Let $E_1$ and $E_2$ be two elliptic curves in the canonical forms defined over $k$

$$E_1: \quad y^2 = x^3 + ax + b$$
$$E_2: \quad y^2 = x^3 + a'x + b'$$

With nonzero discriminates $\Delta \neq 0$.

Let $\phi$ be an isogeny between $E_1$ and $E_2$ and $k(E_1)$ $k(E_2)$ be function fields of $E_1$ and $E_2$ respectively. The degree of $\phi$, $d(\phi)$, can be defined as the degree of corresponding function field extension $k(E_1)/k(E_2)$.

Let $\psi$ be another isogeny from $E_1$ to $E_2$, it is well known that

$$d(\varphi + \psi) + d(\varphi - \psi) = 2d(\varphi) + 2d(\psi) \qquad (*)$$

The relation (*) indicate that $d(\phi)$ is a quadratic function of $\phi$, from this we can deduce that:

Hasse theorem. Let $E$ be an elliptic curve over a finite field $F_q$. The number $N$ of points of $E$ defined over $F_q$ satisfier

$$|N - (q+1)| \leq 2q^{1/2}.$$

Hesse theorem has an important application in cryptography. Many engineers make use of Hasse theorem but cannot totally understand the proof. So it is important and quite interesting to give an elementary proof. The "elementary" we means the knowledge we needed does not suppress undergraduate ,In fact, Hasse theorem can be proved by many methods such as by using of weil pairs , or by using of "dual isogeny" after we induce it. All these proof at least need some, more or less ,abstract concept of mathematics, such as divisor, or concept of Galois extension, It is Manin,in[4], using an idea of Hasse,give an entitely elementary proof of the theorem,the proof of Manin,had been adopt in Knapp book[3] ,In 1971,H.Zimmer [7]presented a valuation theoretic proof ,In [1],Casseles mentioned (*) without a proof, and then in his book[5],J.S.Milne also list the relation (*),with no a proof. In [6], D.Zagier put forward the problem:

Since $$\sum_{x(\bmod p)} \left( \frac{x^3 - 35x + 98}{p} \right) = 0, \qquad \text{when} \quad \left( \frac{p}{7} \right) = -1 \quad \text{and}$$

$= \pm 2A$, when $\left(\dfrac{p}{7}\right) = 1$, $\quad p = A^2 + 7B^2$

D.Zagier asked if we can give an elementary proof for above formula.

It is well known that (*) can be used to get the sum

$$\sum_{x \pmod p} \left( \dfrac{x^3 - 35x + 98}{p} \right)$$

Hence an elementary proof of (*) seems to give a possible answer to Zagier problem.

The present note is devote to give an elementary proof of (*) and then deduce Hasse theorem. We also prove the classical equation satisfied by $\phi$ when $\phi$ is an endomorphism.

2. Some Lemmas

At first we give following

Definition. Let $Q_1, Q_2, \cdots Q_k$ be polynomials in variable $x$, writing $\deg(Q_j)$ as the degree of $Q_j(x)$, we define

$$H(Q_1, Q_2 \cdots, Q_k) = \max(\deg(Q_1), \deg(Q_1), \cdots, \deg(Q_k)).$$

Further if

$$Q(x) = q_n x^n + q_{n-1} x^{n-1} + \cdots + q_0 \quad \text{Where } q_j \in k, q_n \neq 0.$$

We call $q_n$ the leading coefficients of $Q(x)$.

Lemma 1. Let $A, B, C, D$ be polynomials in undetermined $x$ satisfying $(A, B) = 1 \ (C, D) = 1$, where $(\bullet, \bullet)$ represent the great common divisor of the polynomials. Then for $Q_1 = AC, Q_2 = AD + BC, Q_3 = BD$, we have

$$H(Q_1, Q_2, Q_3) = H(A, B) + H(C, D)$$

Proof. At first we assume that $\deg(AC) \geq \deg(BD)$, we divide into three cases to discuss:

Case 1. $\deg(A) \geq \deg(B)$, $\quad \deg(C) \geq \deg(D)$,

In such case, we have $\deg(AC) \geq \deg(AD)$, $\quad \deg(AC) \geq \deg(BC)$.

Hence the polynomials $AC$ is of largest degree. From the definition we have

$$H(Q_1, Q_2, Q_3) = \deg(AC) = \deg(A) + \deg(C) = H(A, B) + H(C, D)$$

Case 2. $\deg(A) > \deg(B)$, $\quad \deg(C) < \deg(D)$

We easily find that $AD > AC$, $AD > BC$, $AD > CD$ so

$$H(Q_1, Q_2, Q_3) = \deg(AD) = \deg(A) + \deg(D) = H(A,B) + H(C,D)$$

Case 3. $\deg(A) < \deg(B)$, $\deg(C) > \deg(D)$

Since $BC > AD$, $BC > BD$, $BC > AC$

$$H(Q_1, Q_2, Q_3) = \deg(BC) = \deg(B) + \deg(C) = H(A,B) + H(C,D)$$

Collect above results we see that the lemma is valid when $\deg(AC) \geq \deg(BD)$.

When $\deg(AC) < \deg(BD)$, the proof is similar so we omit it.

Lemma 2. Let $P, Q, R, S$ be polynomials satisfying $(P,Q) = 1$, $(R,S) = 1$. If

$$Q_1 = (PR - aQS)^2 - 4b(PS + QR)QS$$
$$Q_2 = 2(PR + aQS)(PS + QR) + 4b(QS)^2$$
$$Q_3 = (PS - QR)^2$$

Where $a, b$ are as in elliptic curves (1).

Then

$$H(Q_1, Q_2, Q_3) = 2H(P,Q) + 2H(R,S)$$

Proof. We divide into some cases to prove the Lemma.

(I). $\deg(R) > \deg(S)$. In this case

If $\deg(P) \geq \deg(Q)$, then $\deg(PR) > \deg(QS)$

$$\deg(PR) \geq \deg(QR) \quad \text{And} \quad \deg(PR) \geq \deg(PS)$$

So we get from the definition of $Q_1, Q_2, Q_3$

$$\deg(Q_1) = \deg((PR)^2) = 2\deg(P) + 2\deg(R) = 2H(P,Q) + 2H(R,S)$$

$$\deg(Q_2) \leq \deg((PR)^2) \leq \deg(Q_1)$$

$$\deg(Q_3) \leq \deg((PR)^2) \leq \deg(Q_1)$$

Thus $H(Q_1, Q_2, Q_3) = H(Q_1) = 2H(P,Q) + 2H(R,S)$

The Lemma follows.

If $\deg(P) < \deg(Q)$, similarly we get

$$\deg(Q_3) = \deg((QR)^2) = 2\deg(Q) + 2\deg(R) = 2H(P,Q) + 2H(R,S)$$

$$\deg(Q_1) \leq \deg((QR)^2) \leq \deg(Q_3)$$

$$\deg(Q_2) \leq \deg((QR)^2) \leq \deg(Q_3)$$

Thus $$H(Q_1,Q_2,Q_3) = H(Q_3) = 2H(P,Q) + 2H(R,S)$$

Since $\deg(P) > \deg(Q)$ is symmetric to $\deg(R) > \deg(S)$, the proof is similar we omit, thus we are left with

$\deg(P) \leq \deg(Q)$ and $\deg(R) \leq \deg(S)$

We divide into four cases

(i) $\deg(P) = \deg(Q)$    $\deg(R) < \deg(S)$

(ii) $\deg(R) = \deg(S)$    $\deg(P) < \deg(Q)$

(iii) $\deg(P) < \deg(Q)$    $\deg(R) < \deg(S)$

(iv) $\deg(P) = \deg(Q)$    $\deg(R) = \deg(S)$.

In (i),    we have

$\deg(PS) > \deg(QR)$    Hence

$$\deg(Q_3) = \deg((QS)^2) = 2\deg(Q) + 2\deg(S) = 2H(P,Q) + 2H(R,S)$$

$$\deg(Q_2) \leq \deg((QS)^2) \leq \deg(Q_3)$$

$$\deg(Q_1) \leq \deg((QS)^2) \leq \deg(Q_3)$$

So $$H(Q_1,Q_2,Q_3) = H(Q_3) = 2H(P,Q) + 2H(R,S)$$

Since (ii) is symmetric to (i), we omit its proof.
For case (iii), we have

$\deg(QS) > \deg(PR)$        $\deg(QS) > \deg(PS)$

Let the leading coefficients of $Q, S$ are $l(Q), l(S)$ respectively, then

The leading coefficients of $Q_1$ is $a'^2 l(Q)^2 l(S)^2$,

The leading coefficients of $Q_2$ is $4b' l(Q)^2 l(S)^2$

$\deg(Q_3) \leq 2\deg(QS)$, hence unless $a' = b' = 0$,

$$H(Q_1,Q_2,Q_3) = \deg((QS)^2) = 2\deg(Q) + 2\deg(S) = 2H(P,Q) + 2H(R,S).$$

But  $\Delta = 4a'^3 + 27b'^2 \neq 0$, so  $a' = b' = 0$  is impossible.

In case (iv), Put  $P_1 = PR, P_2 = PS + QR, P_3 = QS$,

Then  $\deg(P_1) = \deg(P_3)$   $\deg(P_2) \leq \deg(P_1)$

We denote the leading coefficients of  $P, Q, R, S$  as  $l(P), l(Q), l(R), l(S)$  respectively, then it is obvious that the leading coefficients of  $P_1, P_2, P_3$  are

$$l(P_1) = l(P)l(R), l(P_3) = l(Q)l(S), l(P_2) = l(Q)l(S) + l(P)l(R)$$

respectively,   furthermore

$$l(P_1) = l(P)l(R) \neq 0, l(P_3) = l(Q)l(S) \neq 0$$

Note that  $(PS - QR)^2 = P_2^2 - 4P_1 P_3$, we rewrite

$$Q_1 = (P_1 - a'P_3)^2 - 4b'P_2 P_3$$
$$Q_2 = 2(P_1 + a'P_3)P_2 + 4b'P_3^2$$
$$Q_3 = P_2^2 - 4P_1 P_3$$

Hence

$$l(Q_1) = (l(P_1) - a'l(P_3))^2 - 4b'l(P_2)l(P_3)$$
$$l(Q_2) = 2(l(P_1) + a'l(P_3))l(P_2) + 4b'(l(P_3))^2$$
$$l(Q_3) = l(P_2)^2 - 4l(P_1)l(P_3)$$

Since  $\deg(Q_j) \leq 2\deg(P_1) = 2\deg(P) + 2\deg(R)$   $j = 1, 2, 3$, Obviously if at least one of  $l(Q_1), l(Q_2), l(Q_3)$  is not zero, then

$$H(Q_1, Q_2, Q_3) = 2\deg(P) + 2\deg(R) = 2H(P, Q) + 2H(R, S)$$

So we suppose that  $l(Q_1) = 0, l(Q_2) = 0, l(Q_3) = 0$.  noting  that  $l(P_3) \neq 0$, Putting

$$x = l(P_2)/l(P_3), \quad y = l(P_1)/l(P_3),$$

From  $l(Q_1) = 0, l(Q_2) = 0, l(Q_3) = 0$  we get

$$y = x^2/4 \triangleq x_1^2$$
$$(y - a')^2 = 4b'x = 8b'x_1$$
$$(y + a')x_1 = -b'$$

Substituting  $y = x^2/4 \triangleq x_1^2$  into other two relations,

we get that $x_1^3 + a'x_1 + b' = 0$,  $x_1^4 - 2a'x_1^2 - 8b'x_1 + a'^2 = 0$

But it is easy to verify that

$$(3x_1^2 + 4a')(x_1^4 - 2a'x_1^2 - 8b'x_1 + a'^2) - (3x_1^3 - 5a'x_1 - 27b')(x_1^3 + a'x_1 + b') = 4a'^3 + 27b'^2 \neq 0$$

we get a contradiction. Thus

$$H(Q_1, Q_2, Q_3) = 2\deg(P) + 2\deg(R) = 2H(P, Q) + 2H(R, S)$$

Colleting above all the Lemma follows.

According to above two lemmas we have

Corollary.  Let  $A, B, C, D$  be as in Lemma 1. Let  $P, Q, R, S$  be as in Lemma 2 . If

$$AC = (PR - a'QS)^2 - 4b'(PS + QR)QS$$
$$AD + BC = 2(PR + a'QS)(PS + QR) + 4b'(QS)^2$$
$$BD = (PS - QR)^2$$

Then   $H(A, B) + H(C, D) = 2H(P, Q) + 2H(R, S)$

Proof . Obvious.

3. Main results

In this section we will prove our main results. At first we repeat some concepts.

Let
$$E_1: \quad y^2 = x^3 + ax + b$$
$$E_2: \quad y^2 = x + a'x + b'$$

Be two elliptic curves over field  $k$  in canonical forms.

Let   $\varphi, \psi$   be two isogenies from   $E_1$   to  $E_2$, defined over the field  $k$ .

Then we have

Theorem 1 Let  $d(\varphi), d(\psi), d(\varphi + \psi), d(\varphi - \psi)$ denote the degree of isogenies of

$\varphi, \psi, \varphi + \psi, \varphi - \psi$   respectively, then

$$d(\varphi + \psi) + d(\varphi - \psi) = 2d(\varphi) + 2d(\psi)$$

Proof. The process is classical.

Put  $(x_1, y_1) = \phi(x, y)$        $(x_2, y_2) = \psi(x, y)$

$(x_3, y_3) = (\phi + \psi)(x, y)$        $(x_4, y_4) = (\phi - \psi)(x, y)$

Thus we get from addition formula of elliptic curve

$$(x_1 - x_2)^2(x_3 + x_4) = 2(x_1 x_2 + a')(x_1 + x_2) + 4b'$$
$$(x_1 - x_2)^2(x_3 x_4) = x_1^2 x_2^2 - 2a'x_1 x_2 - 4b'(x_1 + x_2) + a'^2$$
(**)

We write $x_1 = P(x)/Q(x)$  $x_2 = R(x)/S(x)$, $x_3 = A(x)/B(x)$  $x_4 = C(x)/D(x)$

Where  $(P,Q)=1$  $(R,S)=1$  $(A,B)=1$  $(C,D)=1$

Substituting $x_1$  $x_2, x_3$  $x_4$ into (**) we have

$$\frac{AC}{BD} = \frac{(PR-a'QS)^2 - 4b'(PS+QR)QS}{(PS-QR)^2}$$

$$(AD+BC)/(BD) = (2(PR+a'QS)(PS+QR) + 4b'(QS)^2)/(PS-QR)^2$$

Without loss of generality we put

$$AC = U((PR-a'QS)^2 - 4b'(PS+QR)QS)$$
$$BD = U(PS-QR)^2$$

Where $U$ is a polynomials defined over $k$, Let $F$ denote the prime factor of $U$, then

$$F \mid AC, \quad F \mid BD$$

Assume that  $F \mid A$,  which $F$ cannot divide $B$, so we get $F$ divide $D$. Since

$F \mid AD+BC$  We have  $F \mid BC \Rightarrow F \mid C$,  which is  $F \mid (C,D)=1$.

So we see that $U$ is a constant. hence without of loss of generality we suppose that $U=1$. thus we have that

$$AC = (PR-a'QS)^2 - 4b'(PS+QR)QS$$
$$AD+BC = 2(PR+a'QS)(PS+QR) + 4b'(QS)^2$$
$$BD = (PS-QR)^2$$

From the corollary we find that

$$d(\varphi+\psi) + d(\varphi-\psi) = 2d(\varphi) + 2d(\psi)$$

Corollary 2. A function $f: M \to K$ from a group into field $K$ of characteristic $\neq 2$, satisfying

$$f(x+y) + f(x-y) = 2f(x) + 2f(y)$$

Define $L(x,y) = f(x+y) - f(x) - f(y)$. Then the function $L(x,y)$ is bilinear.

Proof. See, for example, [1] or [5].

4. Hasse theorem and others

In this section we assume two elliptic curves $E_1$, $E_2$ are the same. Let

$$\varphi, \psi: \quad E \to E$$

Both are isogenies. By abuse of notation we denote the constant endomorphism by 0 and the identify endomorphism by 1.

Then from Corollary 2 it is easy to see that

$$d(m\varphi + n\psi) = m^2 d(\varphi) + n^2 d(\psi) + mnL(\varphi,\psi) \qquad (4.1)$$

Where $\quad L(\varphi,\psi) = d(\varphi+\psi) - d(\varphi) - d(\psi)$

Since $\quad d(m\varphi + n\psi) \geq 0 \quad$ then

$$L(\varphi,\psi)^2 \leq 4d(\varphi)d(\psi). \qquad (4.2)$$

Lemma 3. Every endomorphism $\varphi$ satisfies a quadratic equation

$$\varphi^2 - L(\varphi,1)\varphi + d(\varphi) = 0.$$

Proof. We give a proof by elementary method. This is almost the same as [1]. Let $\varphi' = L(1,\varphi) - \varphi$, since $L(1,\varphi) \in Z$ we see that $\phi'$ is also an endomorphism.

From (4.1) we get

$$d(\varphi') = d(L(1,\varphi) - \varphi) = L(1,\varphi)^2 - L(1,\varphi)L(1,\varphi) + d(\varphi) = d(\varphi)$$

Also $d(1+\varphi') = 1 + L(1,\varphi') + d(\varphi')$

$d(1+\varphi) = 1 + L(1,\varphi) + d(\varphi)$. Noting that

$$L(1,\varphi) + L(1,\varphi') = L(1,\varphi+\varphi') = (1+L(1,\varphi))^2 - L(1,\varphi)^2 - 1 = 2L(1,\varphi)$$

So $L(1,\varphi) = L(1,\varphi')$

We have $d(1+\varphi) = d(1+\varphi')$

And $\quad d((1+\varphi)(1+\varphi')) = (1 + L(1,\varphi) + d(\varphi))^2 \qquad (4.3)$

But $d((1+\varphi)(1+\varphi') = d(1+\varphi+\varphi'+\varphi\varphi') = d(1+L(1,\varphi)+\varphi\varphi')$,

Let $m = 1 + L(1,\varphi)$, then

$$d(1+L(1,\varphi)+\varphi\varphi') = d(m+\varphi\varphi') = m^2 + mL(1,\varphi\varphi') + d(\varphi\varphi')$$

Compare with (4.3) we get $L(1,\varphi\varphi') = 2d(\varphi)$. So for any rational integer $l$

$$d(l - \varphi\varphi') = l^2 - lL(1,\varphi\varphi') + d(\varphi)^2 = (l - d(\varphi))^2.$$

If $l = d(\varphi)$, we get the Lemma.

Especially for Frobenius map $\pi(x, y) = (x^q, y^q)$, we have

$\pi^2 - L(\pi,1)\varphi + q = 0$. Where $L(-\pi,1) = d(1-\pi) - q - 1$.

By(4.2) we get

$|d(1-\pi) - q - 1| \leq 2\sqrt{q}$.

Hence the Hasse theorem follows.